\documentclass[11 pt, reqno]{article}

\usepackage{amsmath,amsfonts,amssymb,amsthm}
\usepackage[colorlinks=true,hyperindex=true]{hyperref}
\usepackage{cancel,bbm}
\usepackage{mathabx} 
\usepackage{comment}
\usepackage{cite}
\usepackage{enumerate}

\usepackage{chngcntr}
\counterwithin*{equation}{section}

\usepackage[a4paper, left=2cm, right=2cm, top=2 cm, bottom= 2 cm]{geometry}
\usepackage{color}
\usepackage{fancyhdr}
\usepackage{latexsym}

\newtheorem{ccounter}{ccounter}[section]
\newtheorem{thm}[ccounter]{Theorem}
\newtheorem{lem}[ccounter]{Lemma}
\newtheorem{cor}[ccounter]{Corollary}
\newtheorem{defn}[ccounter]{Definition}
\newtheorem{prop}[ccounter]{Proposition}
\newtheorem{ass}[ccounter]{Assumption}
\newtheorem{ex}[ccounter]{Example}

\def\bet{\begin{thm}}
\def\eet{\end{thm}}
\def\bel{\begin{lem}}
\def\eel{\end{lem}}
\def\bas{\begin{ass}}
\def\eas{\end{ass}}
\def\bec{\begin{cor}}
\def\eec{\end{cor}}
\def\bed{\begin{defn}}
\def\eed{\end{defn}}
\def\bep{\begin{prop}}
\def\eep{\end{prop}}
\def\beq{\begin{equation}}
\def\eeq{\end{equation}}
\def\proof{\noindent {\bf Proof.}\ \ }
\def\bea{\begin{equation*}}
\def\eea{\end{equation*}}

\def\bex{\begin{ex}}
\def\eex{\end{ex}}

\def\rr{\mathbb{R}}

\def\1{\boldsymbol{1}}
\def\Im{\mathrm{Im}}
\def\Re{\mathrm{Re}}
\def\e{\mathrm{e}}
\def\i{\mathrm{i}}
\def\del{\partial}
\def\d{\mathrm{d}}
\def\eps{\varepsilon}
\renewcommand\leq\varleq
\renewcommand\geq\vargeq
\def\ee{\mathrm{E}}

\def\O{\mathcal{O}}

\def\ee{\mathbb{E}}

\def\pp{\mathbb{P}}

\def\GKA{\mathcal{G}_{K, A}}
\newcommand{\Sb}{\mathcal{S}_b}
\def\Fxi{\mathcal{F}_{\xi}}
\def\JD{\mathcal{J}_D}

\def\A{\mathcal{A}}
\def\rhosc{\rho_\mathrm{sc}}
\def\msc{m_{\mathrm{sc}}}
\def\ss{\mathbb{S}}
\def\hatgam{\hat{\gamma}}

\def\tilgam{\tilde{\gamma}}
\def\tilgams{\tilde{\gamma}_s}
\def\Gama{\Gamma_a}
\def\mfb{\mathfrak{b}}
\def\Var{\mathrm{Var}}

\begin{document}
\title{Deformed GOE}


\begin{table}
\centering

\begin{tabular}{c}
\multicolumn{1}{c}{\Large{\bf Free energy fluctuations of the $2$-spin spherical SK model}}\\
\multicolumn{1}{c}{\Large{\bf at critical temperature}}\\
\\
\\
\end{tabular}
\begin{tabular}{c }
Benjamin Landon\\
\\
 Department of Mathematics \\
 \small{MIT} \\
\small{blandon@mit.edu}  \\
\\
\end{tabular}
\\
\begin{tabular}{c}
\multicolumn{1}{c}{\today}\\
\\
\end{tabular}

\begin{tabular}{p{15 cm}}
\small{{\bf Abstract:}  We investigate the fluctuations of the free energy of the $2$-spin spherical Sherrington-Kirkpatrick model  at critical temperature $\beta_c = 1$.  When $\beta = 1$ we find asymptotic Gaussian fluctuations with variance $\frac{1}{6N^2} \log(N)$, confirming in the spherical case a physics prediction for the SK model with Ising spins.  We furthermore prove the existence of a  critical window on the scale $\beta = 1 +\alpha \sqrt{ \log(N) } N^{-1/3}$.  For any $\alpha \in \rr$ we show that the fluctuations are at most order $\sqrt{ \log(N) } / N$, in the sense of tightness.  If $ \alpha \to \infty$ at any rate as $N \to \infty$ then, properly normalized, the fluctuations converge to the Tracy-Widom$_1$ distribution.  If $ \alpha \to 0$ at any rate as $N \to \infty$  or $ \alpha <0$ is fixed, the fluctuations are asymptotically Gaussian as in the $\alpha=0$ case.   In determining the fluctuations, we apply a recent result of  Lambert and Paquette \cite{lambert2020strong} on the behavior of the Gaussian-$\beta$-ensemble at the spectral edge.}
\end{tabular}
\end{table}

\section{Introduction and main results}

The Hamiltonian of the $2$-spin spherical Sherrington-Kirkpatrick (SSK) model is given by,
\beq
H_N ( \sigma) =   \frac{1}{ \sqrt{2N}} \sum_{i , j } \sigma_i g_{ij} \sigma_j 
\eeq
where the $\{ g_{ij} \}_{i, j}$ are independent standard normal random variables and the vector of spins lies in the state space $\sigma \in \ss^{N-1} := \{ \sigma \in \rr^N : ||\sigma||_2^2 = N \}$.  The SSK model was introduced by Kosterlitz, Thouless and Jones \cite{KTJ} as a simplification of the usual SK model which has the same form except that the spins $\sigma$ are assumed to lie in the hypercube.  The model  with Ising spins  was introduced in 1975 in order to explain strange phenomena of various alloys \cite{sherringtonkirkpatrick}.  Since its introduction, this model and its generalizations have been the subject of much research in both the physics and mathematics communities.  We refer the interested reader to, e.g., \cite{mezard1987spin,tal1,tal2,panchenko2013sherrington}.

One of the primary thermodynamic quantities of interest is the free energy $F_N ( \beta)$ which is defined by,
\beq
F_N ( \beta ) := \frac{1}{N} \log (Z_N ( \beta )) , \qquad Z_N (\beta ) := \frac{1}{ | \ss^{N-1} |} \int_{\ss^{N-1}} \e^{ - \beta H_N ( \sigma) } \d \omega_{N-1} ( \sigma)
\eeq
where $\d \omega_{N-1}$ is the uniform surface measure on $\ss^{N-1}$.  Here, $Z_N (\beta)$ is the partition function, the normalizing constant of the Gibbs measure on $\ss^{N-1}$ that has weight proportional to $\e^{ - \beta H_N ( \sigma)}$.  Both the SK and SSK models are known to exhibit a phase transition at the critical temperature $\beta_c :=1$.  In the case of the SSK, the limit of the free energy was calculated by \cite{CS} to be, 
\beq
\lim_{N \to \infty} F_N ( \beta ) = f ( \beta ) := \begin{cases} \beta^2/4 , & \beta \leq 1 \\ \beta - \frac{\log ( \beta)}{2} - \frac{3}{4} & \beta \geq 1 \end{cases}.
\eeq
This was later rigorously proven by Talagrand \cite{talagrand-FE}.

 Baik and Lee \cite{baiklee} further described the fluctuations of the free energy around its limiting value for non-critical $\beta \neq 1$ as follows.  In the high temperature regime $\beta < 1$ the fluctuations are $\O(N^{-1})$ and are asymptotically Gaussian.  In the low temperature regime $\beta >1$, the fluctuations are $\O (N^{-2/3})$ and converge to the Tracy-Widom$_1$ distribution (TW$_1$) of random matrix theory, the limiting distribution of the largest eigenvalue of a symmetric  random matrix with Gaussian entries \cite{tracywidom} (note that \cite{baiklee} uses a different scaling in which the transition between high and low temperature regimes is at $\frac{1}{2}$).  The high temperature central limit theorem for the SK model with Ising spins was proven earlier by Aizenman, Lebowitz and Ruelle \cite{aizenmanlebowitzruelle}.  

Less is known about the size or nature of the fluctuations of either the SK or SSK models at or near the critical temperature $\beta_c = 1$.   Talagrand examined the SK model on the scale $\beta = 1  +c N^{-1/3}$ and found that the overlap undergoes a phase transition depending on whether or not $c$ diverges \cite{tal1}.  At $\beta_c$ he showed that the variance is at most $\O (N^{-3/2})$.  Chatterjee proved a general estimate that for fixed $\beta$ the variance of the free energy of the SK model is at most $\O ((N \log(N) )^{-1})$ \cite{chatterjee2009disorder}.  Chen and Lam \cite{chen2019order} recently obtained the estimate for the SK model,
\beq
\Var ( F^{(\mathrm{SK})}_N ( \beta_c ) ) \leq C \frac{ \log(N)^2}{ N^2},
\eeq
which is the tightest bound for the variance at critical temperature.

For the SK model, it is predicted that \cite{aspelmeier2008free,parisi2009phase}
\beq \label{eqn:sk-pred}
\Var (  F^{(\mathrm{SK})}_N ( \beta_ c) ) = \frac{ \log(N) }{6N^2} + \O (1).
\eeq
Our first main result addresses this prediction for the SSK model.
\bet \label{thm:crit-temp}
Let $\beta = \beta_c=1$ and $F_N(\beta)$ be the free energy of the SSK model, with $f( \beta)$ its limiting value as above. Then, the random variable
\begin{align}
\frac{ N(F_N ( \beta_c ) - f( \beta_c)) + \frac{1}{12} \log(N) }{\sqrt{ \frac{1}{6} \log(N) }}
\end{align}
converges in distribution to a standard normal random variable as $N \to \infty$.
\eet
Based on the form of the variance that they find in the high and low temperature regimes, Baik and Lee \cite{baiklee} predicted that there is a critical regime in the SSK model near the critical temperature of the form
\beq
\beta = 1 + \alpha \sqrt{ \log(N)} N^{-1/3} 
\eeq
for $\alpha \in \rr$ of order $1$, for which the size of the fluctuations are  $\O ( N^{-1} \sqrt{ \log(N) })$.   Our next theorem confirms this prediction by showing that the fluctuations are at most of this order.  Moreover, we find that in the cases $(\alpha)_+ \to 0$ and $\alpha \to \infty$, the free energy fluctuations converge to the Gaussian and Tracy-Widom$_1$ distributions that characterize the high and low temperature regimes, respectively. 
\bet \label{thm:main}
Let $\beta = 1 + \alpha \sqrt{ \log(N) } / N^{1/3}$.  Then the random variable,
\beq
\frac{ N(F_N ( \beta ) - f( \beta)) + \frac{1}{12} \log(N) }{\sqrt{ \frac{1}{6} \log(N) }}
\eeq
is tight.  If $\alpha \to 0$ or $ \alpha <0$ is fixed, then it converges to a standard normal random variable. Suppose that $\alpha \to \infty$ as $N\to \infty$ but slowly enough that $\beta$ remains bounded above. Then,
\beq
\frac{1}{N^{1/3} (\beta -1 ) } \left( NF_N ( \beta) - f ( \beta)) +  \log(N)/12 \right)
\eeq
converges to a Tracy-Widom random variable.  
\eet
Note that Theorem \ref{thm:crit-temp} follows from the above result.  The restriction that $\beta$ is bounded above is  for convenience; it is likely that the result holds without this condition under minimal modification of our methods or those of Baik-Lee \cite{baiklee}. We are interested in the behavior of $\beta$ near $1$ so we do not pursue this minor extension in our work.  

 In many respects,  the behavior of the SSK model turns out to be quite different from the SK model with Ising spins or  higher-order spin models with either continuous or discrete state spaces.  This is primarily due to the fact that the combination of the quadratic nature of the SSK model and the continuous nature of the state space reduces the complexity of the {\it energy landscape}, compared to other ``true'' spin glasses.  For more complex models the number of critical points of the Hamiltonian is typically exponential in $N$ whereas for the SSK it is linear - see, e.g., \cite{auffinger2013random}.  Moreover, the SSK model is usually studied using random matrix methods instead of approaches more commonly associated with the theory of spin glasses.

The main simplification arising in the SSK model is a contour integral formula for observables in terms of the eigenvalues of a matrix formed from the disorder random variables,
\beq
H_{ij} = -\frac{ g_{ij} +g_{ji} }{ \sqrt{2 N } }.
\eeq
The matrix $H$ is an element of the Gaussian Orthogonal Ensemble (GOE), one of the main objects of random matrix theory, and accordingly, much is known about its spectral properties.  The contour integral was used by Kosterlitz-Thouless-Jones \cite{KTJ} in their original paper to study the free energy.  Baik and Lee's work \cite{baiklee} further examined this formula and used it to expand the free energy in terms of the eigenvalues of $H$.  In the high temperature regime, the leading order fluctuations of the free energy are given by a linear spectral statistic,
\beq \label{eqn:lss}
\sum_i \varphi ( \lambda_i )
\eeq
for a certain test function $\varphi$.  Such objects are well studied in random matrix theory, see e.g., \cite{lytovapastur,shcherbina2011central} where it is proven for sufficiently regular $\varphi$ that the limiting fluctuations are Gaussian.  In the low temperature case, Baik and Lee showed that the free energy fluctuations are determined by the largest eigenvalue of $H$, whose convergence was established by Tracy and Widom \cite{tracywidom}; the limiting distribution is now called the Tracy-Widom$_1$ distribution. 

Baik-Lee and Baik-Lee-Wu have since extended their analysis of the free energy to other models related to the SSK \cite{baiklee1,baiklee2,baiklee3}, including a bipartite model and models incorporating a Curie-Weiss-type interaction.   Other recent developments have studied the \emph{overlap} between two samples (\emph{replicas}) from the random Gibbs measure as well as the SSK model with a magnetic field.  The overlap is an important quantity in the theory of spin glasses, begin linked to the free energy via the Parisi formula \cite{parisi1,parisi2,parisi3}. In the work \cite{vl-ps}, Nguyen and Sosoe extended the contour integral formula of \cite{KTJ} to prove a central limit theorem for the overlap in the high-temperature regime. 
  Subsequently, the author with Sosoe \cite{SSK-ls} examined the overlap fluctuations in the low temperature regime and found that they are determined by an explicit function of the Airy$_1$ random point field, which is the joint limit of the largest eigenvalues of the GOE.  In the work \cite{landon2020fluctuations} we analyzed the fluctuations of the free energy and overlap for the SSK model with magnetic field in various scaling regimes, confirming some predictions of Fyodorov and Le Doussal \cite{fyodorov2014topology}.  This model was also investigated by Baik, Collins-Wildman, Le Doussal and Wu \cite{baik-wu}.

We now recall the notation $\beta = 1 + \alpha \sqrt{ \log(N) } N^{-1/3}$.   
When $(\alpha)_+ \to 0$ we show that the fluctuations of the free energy are determined by the logarithm of the characteristic polynomial of the matrix $H$ evaluated at a point at the edge of the limiting spectral measure  (in this case at $E = 2 + \O (N^{-2/3})$).  Formally, this has the form of \eqref{eqn:lss} but with a singular function $\varphi$, for which the existing central limit theorems of random matrix theory do not apply.  However, recently Lambert and Paquette \cite{lambert2020strong} have analyzed precisely the quantity we find.  In particular, they prove a central limit theorem for the log-characteristic polynomial, allowing us to obtain Gaussian fluctuations.  When $\alpha >0$ we find that the free energy fluctuations are determined by a sum of the log-characteristic polynomial and the largest eigenvalue of $H$, and so we obtain tightness of the free energy fluctuations.$^1${\let\thefootnote\relax\footnotetext{1. In private communication, Lambert and Paquette indicate that it follows from their methods \cite{lambert2020strong} that the log-characteristic polynomial and the largest eigenvalues of $H$ are asymptotically independent.  Conditional on this, the $\alpha >0$ regime of our main theorem is restated as convergence to an independent sum of Tracy-Widom$_1$ and Gaussian random variables. }}  

The main technical element of the present work is then to relate the free energy to the eigenvalues of $H$ using the method of steepest descent.  In particular, the method of steepest descent has not been carried out for the SSK model at critical temperature.  The main technical complication is that the saddle itself fluctuates on the same scale as the eigenvalues.  In order to handle this we require as input estimates for the eigenvalue positions on their natural scale (estimates in particular stronger than the well-known rigidity estimates and local semicircle law of random matrix theory -see, e.g.,  \cite{KBG}). The  estimates we require were proven with the author with Sosoe in our work \cite{SSK-ls}.  They are proven using a result of Gustavson \cite{gust} on the eigenvalues of the Gaussian Unitary Ensemble (GUE) and a coupling between the GOE and GUE due to Forrester-Rains \cite{FR}.

\vspace{5 pt} 

\noindent{\bf Applications to statistics.}  Consider  the measure $P_\lambda$ on $N \times N$ symmetric matrices induced by
\beq
H_\lambda := H + \lambda v v^T
\eeq
where $H$ is a matrix distributed according to the GOE and $v$ is a vector uniformly distributed on the unit sphere.  This is an example of a spiked random matrix, introduced by Johnstone \cite{johnstone2001distribution} as a simple high-dimensional model of the form ``signal+ noise.''   After its introduction, further investigations have shown the existence of a spectral transition at the point $\lambda_c =1 $ where the largest eigenvalue of $H_\lambda$ separates from the spectral bulk $[-2, 2]$ for larger $\lambda$ (the celebrated Baik-Ben Arous-Pech{\'e} transition \cite{baik2005phase} see also, e.g., \cite{baik2006eigenvalues,feral2007largest} for important work).  Since Johnstone's seminal work a large statistical and mathematical literature on spiked models has emerged; we refer  to, e.g., \cite{el2020fundamental,perry2018optimality} and the references therein for further discussion.

In this context, the quantity
\beq
N (F_N ( \beta ) - \beta^2/4)
\eeq
turns out to be the likelihood ratio $ \frac{P_\lambda}{ P_0}$ under the ``null hypothesis'' $\lambda = 0$.  In the somewhat different setting of sample-covariance matrices, Onatski, Moreira and Hallin \cite{stats_ssk} proved a central limit theorem for the likelihood ratio below the spectral transition (which is at a different point for their model but the behavior is analogous to our setting), implying the mutual contiguity of the null and non-null measures.  They moreover proved that above the spectral transition the likelihood ratio converges to $0$ (in fact is exponentially small).  As a corollary, our work extends this result that the likelihood ratio tends to $0$ as $N \to \infty$ near the spectral transition $\lambda = 1$.

\vspace{5 pt}

\noindent{\bf Organization of paper.}  In the next section we introduce some notation as well as collect the results of random matrix theory that we use in our work.  In Section \ref{sec:saddle} we establish preliminary estimates on the saddle that we use in the method steepest descent as well as analyze the leading order term.  In Section \ref{sec:steep} we carry out the method of steepest descent.  In Section \ref{sec:expand} we prove our technical propositions which expand the free energy of the SSK in terms of random matrix quantities.  Theorem \ref{thm:main} is proven in Section \ref{sec:main}.

\vspace{5 pt}

\noindent{\bf Acknowledgements.}  The author thanks P. Sosoe for discussions on the spherical SK model.  The author thanks E. Paquette and G. Lambert for discussions about the results and methods of their work \cite{lambert2020strong}.

\section{Preliminaries}

\subsection{Notation}

We fix some notation.  When we have a complex parameter $z$ we will denote its real and imaginary parts by
\beq
z = E + \i \eta.
\eeq
The Hamiltonian of the SSK model is
\beq
H_N ( \sigma ) = \frac{1}{ \sqrt{2N}} \sum_{i , j } \sigma_i g_{ij} \sigma_j   =: -\sigma^T H \sigma
\eeq
where $g_{ij}$ are independent standard normal random variables and $H$ has matrix elements,
\beq
H_{ij} = - \frac{ g_{ij} +g_{ji}}{ \sqrt{2N}}.
\eeq
Then $H$ has a random matrix from the Gaussian Orthogonal Ensemble and the relation between the Hamiltonian and $H$ is,
\beq
H_N ( \sigma ) =- \frac{1}{2} \sigma^T H \sigma.
\eeq
We will denote the eigenvalues of $H$ in decreasing order $\lambda_1 \geq \lambda_2 \geq \dots \geq \lambda_N$.  Note that they are almost surely distinct.  For convenient reference we recall the definitions of the free energy and the partition function,
\beq
F_N ( \beta) = \frac{1}{N} \log (Z_N ( \beta)), \qquad Z_N ( \beta ) := \frac{1}{ | \ss^{N-1} | } \int_{ \ss^{N-1}} \e^{ - \beta H_N ( \sigma )  } \d \omega_{N-1} ( \sigma)
\eeq
where the $\ss^{N-1}$ is the $N-1$ dimensional sphere of radius $\sqrt{N}$, $\omega_{N-1}$ is the uniform measure on $\ss^{N-1}$ and so
\beq
| \ss^{N-1} | = \int \d \omega_{N-1} ( \sigma ) = \frac{ 2 \pi^{N/2}}{ \Gamma \left( \frac{N}{2} \right) } N^{\frac{N-1}{2}}.
\eeq

\noindent{\bf Assumption.}  \emph{In general we allow $\beta$ to depend on $N$.  However for definiteness throughout the remainder of the paper we will assume that there is a $\mathfrak{m} >0$ so that,}
\beq
\frac{1}{ \mathfrak{m}} \leq \beta \leq \mathfrak{m}
\eeq
\noindent While it is possible to consider the cases of $\beta \to 0$ or $\infty$ using minor modifications of either our methods or those of Baik-Lee \cite{baiklee} we refrain from doing so as we are mainly interested in the behavior near $\beta = \beta_c = 1$.

An important role will be played by the function,
\beq
G(z) := \beta z - \frac{1}{N} \sum_{i=1}^N \log ( z - \lambda_i)
\eeq
where the $\log$ is the principal branch of the logarithm.  Throughout the paper we will denote the saddle $\gamma$ as the unique solution $\gamma > \lambda_1 $ to,
\beq
G' ( \gamma ) = 0.
\eeq
The almost sure limit of the empirical eigenvalue measure of $H$ is given by Wigner's semicircle distribution,
\beq
\rhosc (E) = \frac{1}{ 2 \pi} \sqrt{4 - E^2 } \1_{ \{ |E| \leq 2 \} },
\eeq
which has Stieltjes transform,
\beq
\msc (z) = \int \frac{1}{ x -z } \rhosc (x) \d x = \frac{ - z + \sqrt{z^2 -4 }}{2}.
\eeq
The empirical Stieltjes transform will be denoted by,
\beq \label{eqn:msc-form}
m_N (z) = \frac{1}{N} \sum_{j=1}^N \frac{1}{ \lambda_i - z }.
\eeq
We say that an event $\A$ (really, a sequence of events $\A = \A_N$) holds with overwhelming probability if for every $D>0$ there is a constant $C_D$ on which
\beq
\pp [ \A] \geq 1- C_D N^{-D}.
\eeq

\subsection{Estimates on eigenvalue positions} \label{sec:eig}

The $N$-quantiles or classical eigenvalue locations $\{\gamma_i \}_{i=1}^N$ of the semicircle distribution are defined by,
\beq
\frac{i}{N}  = \int_{\gamma_i}^{2} \rhosc (E) \d E.
\eeq
A straightforward calculation gives for $i \leq N/2$,
\beq \label{eqn:classical-est}
2 - \gamma_i = \left( \frac{3 \pi i}{2N} \right)^{2/3} + \O \left( \frac{i^{4/3}}{N^{4/3} } \right).
\eeq
Note that the error term is less than the rigidity error appearing in the following result $i^{-1/3} N^{-2/3}$ as long as $i \leq N^{2/5}$.

The following result is from, e.g., \cite{KBG}. 
\bet
For any $\xi >0$ define the event  $\Fxi$ by
\beq
\Fxi := \bigcap_{1 \leq i \leq N } \left\{ | \lambda_i - \gamma_i | \leq \frac{N^{\xi}}{N^{2/3} \min\{ i^{1/3}, (N+1-i)^{1/3} \} } \right\}.
\eeq
Then $\Fxi$ holds with overwhelming probability.
\eet

The following is a result of Section 6.2 of \cite{SSK-ls} as well as the fact that $N^{2/3} ( \lambda_1 -2)$ converges to a random variable. 
\bel \label{lem:gka}
For $K>0$ and $A>0$ let us denote the event $\GKA$, by
\beq
\GKA :=   \left\{ \bigcap_{K \leq j \leq N^{2/5}} \left\{   \left| N^{2/3} ( \lambda_j -2 ) + \left( \frac{ 3 \pi j }{2} \right)^{2/3}  \right| \leq \frac{j^{2/3}}{10} \right\} \right\} \bigcap \{ N^{2/3} ( \lambda_1 - 2) \leq A \}.
\eeq
Let $\eps >0$.  Then there are $K, A$ so that  $\pp[ \GKA] \geq 1- \eps$ for all $N$ large enough.
\eel

Since $\lambda_2$ and $\lambda_1$ converge jointly to the largest particles of the Airy$_1$ random point field which is almost surely simple (see, e.g., Proposition 3.5 of \cite{ramirez2011beta}) we have the following.
\bel \label{lem:sb}
Define for $ b>0$ the event
\beq
\Sb := \{  N^{2/3} ( \lambda_2 - \lambda_1 ) > b \}.
\eeq
Let $\eps >0$.  There is a $b>0$ so that $\pp [ \Sb] \geq 1 - \eps$ for all $N$ large enough.
\eel
We have also the following, proven in Section 6.2 of \cite{SSK-ls}.  
\bel \label{lem:expect-est}
There are positive $C, C_1$ so that,
\beq
\ee\left[ \1_{ \{ N^{2/3} ( \lambda_j -2 ) \leq - C_1 \} } \left| N^{2/3} ( \lambda_j - 2) + \left( \frac{3 \pi k }{2} \right)^{2/3} \right| \right] \leq \frac{C \log(j)^2}{j^{1/3}}, \qquad j \leq N^{2/5}.
\eeq
\eel
The following is a consequence of Theorem 6.1 of \cite{SSK-ls}.
\bel  \label{lem:jd}
For any $D>0$,  let $\JD$ be the event that
\beq
N^{1/3} \left| 1 + \frac{1}{N} \sum_{j=2}^N \frac{1}{ \lambda_j - \lambda_1 } \right| \leq D
\eeq
For any $\eps >0$ there is a $D>0$ so that $\pp [ \JD ] \geq 1- \eps$ for all $N$ large enough.
\eel

We now prove the following simple consequence of the above lemmas.

\bel \label{lem:sad-1}
Let $A, K>C$ for $C$ sufficiently large.   There is a $C'> 0$ so that for $N^{2/3}(E - 2 ) \geq 2 A $ we have for all small $1/100 \geq \xi >0$,
\begin{align} \label{eqn:rig-est}
\ee\left[ \1_{\GKA \cap \Fxi } N^{1/3} \left| \frac{1}{N}\sum_{j=1}^N \frac{1}{ \lambda_j - E } - \msc (E) \right| \right] &\leq \frac{ C' (K + \log^2 (N^{2/3} (E-2) )}{ N^{2/3} (E-2 ) } \notag\\
&+ C' \left(N^{-2/3+\xi} + \1_{ \{ N^{2/3}(E-2 ) \geq N^{1/15-\xi} \}} \frac{N^{\xi}}{ N^{2/3} (E-2) } \right)
\end{align}
\eel
\proof We write,
\begin{align}
 \frac{1}{N}\sum_{j=1}^N \frac{1}{ \lambda_j - E } - \msc (E) &= \sum_{i>K}^N \int_{\gamma_{i-1}}^{ \gamma_i } \frac{ (x - \gamma_i ) + ( \gamma_i - \lambda_i ) }{ ( \lambda_i - E )( x - E) } \d E \notag \\
&+\frac{1}{N} \sum_{i=1}^{K} \frac{1}{ \lambda_i - E} + \int_{0}^{\gamma_{K}} \frac{ \rhosc (x) }{ x- E} \d x.
\end{align}
As long as $N^{2/3} (E  -2 ) > 2 A$, the denominators in the last line are bounded below by $N^{2/3}(E-A)/2$  on the event $\GKA$.  Therefore,
\beq
\1_{\GKA} N^{1/3} \left| \frac{1}{N} \sum_{i=1}^{K} \frac{1}{ \lambda_i - E} + \int_{0}^{\gamma_{K}} \frac{ \rhosc (x) }{ x- E} \d x \right| \leq \frac{CK}{N^{2/3} (E - A)}.
\eeq
As long as $K$ is sufficiently large, then $N^{2/3} ( \lambda_j - 2 ) \geq  c j^{2/3}$ on $\GKA \cap \Fxi$, for all $j \geq K$.  From this we see that,
\begin{align}
& N^{1/3} \1_{ \GKA \cap \Fxi} \left| \sum_{i=K}^N \int_{\gamma_{i-1}}^{ \gamma_i } \frac{ (x - \gamma_i ) + ( \gamma_i - \lambda_i ) }{ ( \lambda_i - E )( x - E) } \d E \right| \notag\\
 \leq & C \1_{ \GKA}  \sum_{j = K}^{N^{1/10}}  \left( j^{-1/3} + \left| N^{2/3} ( \lambda_j - 2) + \left( \frac{ 3 \pi j}{2} \right)^{2/3} \right| \right) \frac{1}{ j^{4/3} + ( N^{2/3} (E -2 ))^2 } \notag\\
+ & C \sum_{ j > N^{1/10}}  \frac{1}{ j^{4/3} +  ( N^{2/3} (E-2) )^2 } \frac{N^{\xi}}{\min\{ j^{1/3}, (N+1-j)^{1/3} \} }
\end{align}
Note we used \eqref{eqn:classical-est}.  For the last line, the contributions for terms in the sum for $j > N/2$ are $\O(N^{-2/3+\xi})$. The remaining terms contribute $\O(N^{-15+\xi})$ if $N^{2/3}(E-2) \leq N^{1/15}$ which can be absorbed into the first line of \eqref{eqn:rig-est} if $N^{2/3}(E-2) \leq N^{1/15-\xi}$ and the second otherwise.  If $N^{2/3}(E-2) \geq N^{1/15}$ then the sum can be further estimated by dividing into cases when $j^{2/3}$ is larger or less than $N^{2/3}(E-2)$.  Ultimately, this is bounded by the second line of \eqref{eqn:rig-est}.

Taking $K$ large enough so that on $\GKA$ we have that $N^{2/3} ( \lambda_j -2 ) \leq -C_1$ where $C_1$ is from Lemma \ref{lem:expect-est}, we see
\begin{align}
& \ee\left[ \1_{ \GKA}  \sum_{j = K}^{N^{1/10}}  \left( j^{-1/3} + \left| N^{2/3} ( \lambda_j - 2) + \left( \frac{ 3 \pi j}{2} \right)^{2/3} \right| \right) \frac{1}{ j^{4/3} + ( N^{2/3} (E -2 ))^2 } \right] \notag\\ 
\leq& C \sum_{j=K}^{N^{1/10}} \frac{ j^{-1/3} \log(j)^2}{ j^{4/3} + (N^{2/3} (E-2))^2 } \notag\\
\leq & C \sum_{j^{2/3} < N^{2/3} (E-2) } \frac{ j^{-1/3} \log(j)^2}{ j^{4/3} + (N^{2/3} (E-2))^2 }  + C \sum_{j^{2/3} \geq N^{2/3} (E-2) } \frac{ j^{-1/3} \log(j)^2}{ j^{4/3} + (N^{2/3} (E-2))^2 }  \notag\\
\leq & C \frac{ \log^2 (N^{2/3} ( E - 2 ) )}{ N^{2/3} (E-2 ) } 
\end{align}
 This concludes the claim. \qed

We require also the following.
\bel \label{lem:msc-behaviour}
For $10 \geq E \geq 0$ and $0 \leq \eta \leq 10$ we have,
\beq
c\sqrt{ |E-2| + \eta } \leq \left| 1 + \msc (z) \right| \leq C \sqrt{ |E-2| + \eta }.
\eeq
For $2 \leq E \leq 10$ we have,
\beq
- \msc \leq 1 - c \sqrt{E-2}.
\eeq
\eel
\proof The first follows from the estimate (4.2) of \cite{erdHos2013local} if we can show that $c \leq |-1-\msc(z)| \leq C$ for $z$ as specified above.  This is easily verified using \eqref{eqn:msc-form}.  The second follows from the first estimate together with the fact that $\msc(2) = -1$ and that $\msc'(E) >0$ for $E > 2$. \qed

\subsection{Convergence result}

For any $Q>0$ define the random variable
\begin{align}
X_Q := \sum_{i=1}^N \log |2 + QN^{-2/3} - \lambda_i | - \frac{N}{2} - N^{1/3} Q + \frac{1}{6} \log(N).
\end{align}
The following is Corollary 1.2 of \cite{lambert2020strong}.  Note that they use a different scaling, where the eigenvalues of the GOE are asymptotically in the interval $[-1, 1]$. 
\bet \label{thm:pq}
For any $Q>0$ we have that
\beq
\frac{X_Q}{ \sqrt{2 \log(N)/3}} 
\eeq
converges to  a standard normal random variable.
\eet

\subsection{Integral representation}

Due to \cite{baiklee} we have for any matrix $H$,
\beq
\int \exp \left[ \frac{ \beta}{2} \sigma^T H \sigma \right] \d \omega_{N-1} ( \sigma ) = \frac{ \beta N^{1/2}}{ 2 \pi \i } \left( \frac{ 2 \pi }{ \beta } \right)^{ \frac{N}{2}} \int_{ a- \i \infty} ^{a + \i \infty} \exp \left[ \frac{N}{2} G(z) \right] \d z
\eeq
for any $a > \lambda_1 (H)$.  Hence we have,
\bel \label{lem:rep}
The following representation for the free energy holds.
\begin{align}
F_N ( \beta ) &= \frac{1}{N} \log (Z_N ( \beta ) ) = \frac{1}{N} \log \Gamma \left( \frac{N}{2} \right)  - \frac{N-2}{2N} \log ( \beta N ) + (\frac{N}{2} - 1) \log(2)  \notag\\
&+ \log \left( \frac{1}{ 2 \pi \i } \int_{a - \i \infty}^{a + \i \infty} \exp \left[ \frac{N}{2} G (z) \right] \d z\right)
\end{align}
where $a > \lambda_1$.
\eel

Note that due to Stirling's formula, for $|1-\beta| \leq 1/2$, 
\begin{align} \label{eqn:stir}
\log \Gamma (N/2) - \frac{N-2}{2} \log (\beta N) = \frac{1}{2} \log(N) - \frac{N}{2} \log (2) - \frac{N}{2} \log ( \beta) - \frac{N}{2} + \O (1).
\end{align}

\section{Estimates on saddle and $G ( \gamma)$}
\label{sec:saddle}

Throughout this section we will make use of the events $\Sb, \JD, \GKA$ and $\Fxi$ that were defined in Section \ref{sec:eig}. 
Recall that the saddle $\gamma$ the unique solution satisfying $\gamma > \lambda_1$ to the equation
\beq
G' ( \gamma ) = 0.
\eeq
Note that this is the solution to,
\beq
\beta = \frac{1}{N} \sum_{i=1}^N \frac{1}{ \gamma - \lambda_i }.
\eeq
Indeed, note that for $\gamma > \lambda_1$ the RHS is a monotonically decreasing function that goes to $\infty$ as $\gamma \to \lambda_1$ and $0$ as $\gamma \to \infty$, guaranteeing a unique solution. 

We first prove the following general estimate.
\bep \label{prop:sad-est} Assume $\beta > c>0$.  
Let $ \eps >0$.  There is a constant $C_\eps >0$ so that,
\beq
\frac{1}{C_\eps(1 + N^{1/3}(\beta-1)_+)} \leq N^{2/3} ( \gamma - \lambda_1 ) \leq C_\eps (1 +  (N^{1/3} ( 1-\beta )_+)^2 ).
\eeq
with probability at least $1 - \eps$ for all $N$ large enough.
\eep

\bel \label{lem:sad-2}
Let $b, A, K, D>0$, with $K$ sufficiently large.  On the event $\GKA \cap \JD \cap \Sb \cap \Fxi$ with $\xi = 1/100$ we have for $E  > \lambda_1$,
\begin{align}
\frac{1}{N} \sum_{i=1}^N \frac{1}{ E- \lambda_i } \geq 1 - DN^{-1/3}  - C' N^{-1/3} (K b^{-1}    + N^{2/3} (E - \lambda_1 ) ) + \frac{1}{ N^{1/3}} \frac{1}{ N^{2/3} ( \lambda_1 - E ) }.
\end{align}
\eel
\proof  We write,
\begin{align}
1 + \frac{1}{N} \sum_{i=1}^N\frac{1}{\lambda_i - E} &=  \left( 1 + \frac{1}{N} \sum_{i=2}^N \frac{1}{ \lambda_i - \lambda_1 } \right) \notag\\
&+ \left( \frac{1}{N} \sum_{i=2}^N \frac{1}{ \lambda_i - E } - \frac{1}{ \lambda_i - \lambda_1 } \right) \notag\\
& + \frac{1}{N^{1/3}}  \frac{1}{ N^{2/3} ( E- \lambda_1  ) }.
\end{align}
Note that for $K$ large enough we have on $\GKA \cap \Fxi$ (with $\xi = 1/100$) that $N^{2/3} (\lambda_{j} - \lambda_1) >c j^{2/3}$ for $j \geq 10 K$. Hence, the term on the second line is bounded by,
\begin{align}
 \left| \frac{1}{N} \sum_{i=2}^N \frac{1}{ \lambda_i - E } - \frac{1}{ \lambda_i - \lambda_1 } \right| &\leq \frac{C K}{ N^{1/3} b} + N^{2/3}(E-\lambda_1) \sum_{j >1 } \frac{1}{ N^{1/3} j^{4/3}}.
\end{align}
This yields the claim. \qed

\noindent{\bf Proof of Proposition \ref{prop:sad-est}}.   From Lemma \ref{lem:sad-1} and Markov's inequality, we see that for any $\eps >0$, there is a $C_2 >0$  so that for any $N^{2/3} (E - 2 ) \geq C_2$ there is an event (depending on $E$) on which $\lambda_1 \leq E$ and
\beq
\frac{1}{N} \sum_{j=1}^N \frac{1}{ E - \lambda_j } \leq N^{-1/3} - \msc (E) \leq N^{-1/3} +\beta + (1-\beta) - cN^{-1/3} \sqrt{ N^{2/3} (E-2 ) },
\eeq
where we used Lemma \ref{lem:msc-behaviour} in the second estimate. 
For $N^{2/3} (E-2) \geq C (1 + N^{2/3}(1-\beta)_+^2 )$ for some $C>0$, the RHS is less than $\beta$.  This proves the upper bound we require.  The lower bound follows immediately from Lemma \ref{lem:sad-2}. \qed

The following lemma will be useful. 
\bel \label{lem:gpp}
On $\GKA \cap \Fxi$ for $K$ sufficiently large and $\xi < 1/100$ we have for $E > \max\{ 2, \lambda_1\}$ that,
\beq
\frac{1}{N^{4/3}} \sum_j \frac{1}{ ( \lambda_j - E)^2} \leq \frac{K}{ N^{4/3} (E-\lambda_1)^2} + \frac{C}{ \sqrt{N^{2/3} (E-2) }}.
\eeq
If $\Sb$ holds,
\beq
\frac{1}{N^{4/3}} \sum_{j>1} \frac{1}{ ( \lambda_j - \lambda_1 )^2} \leq C_{b, k} 
\eeq
\eel
\proof For $j > K$ it holds that $E- \lambda_j \geq  ( E- 2 ) + ( 2- \lambda_j ) \geq (E- \lambda_1 ) + c N^{-2/3} j^{2/3}$ if $K$ is sufficiently large.  Hence,
\begin{align}
\frac{1}{N^{4/3}} \sum_j \frac{1}{ ( \lambda_j - E)^2}  \leq \frac{K}{ N^{4/3}(E- \lambda_1)^2} + \sum_{j > K} \frac{1}{ j^{4/3} + N^{4/3}(E -2)^{2}}
\end{align}
The second sum is easily seen to be $\O ( (N^{2/3} (E- 2))^{-1/2} )$, by dividing it into the cases $j^{2/3} < N^{2/3} (E-2)$ and $j^{2/3} \geq N^{2/3} (E-2)$.  The second estimate is similar. \qed

We now proceed slightly differently in the high temperature and low temperature cases to find better estimates.
\subsection{High temperature $\beta \leq 1$}
In the high temperature case we define $\tilgam$ by
\beq
\beta + \msc ( \tilgam) = 0.
\eeq
Due to Lemma \ref{lem:msc-behaviour} we have in the case that $\beta <1$ that,
\beq
c (1 -\beta)^2 \leq \tilgam -2 \leq C (1-\beta)^2.
\eeq
We first need the following rough bound.
\bel \label{lem:ht-saddle-bd}
Suppose that $ \beta \leq 1$.  For any $\eps >0$ there is a $C_\eps$ for with probability at least $1-\eps$ for all large $N$ we have te
\beq
(C_\eps)^{-1}(1 + N^{2/3} ( 1 - \beta )^2 ) \leq N^{2/3} ( \gamma -2 ) \leq C_\eps ( 1 + N^{2/3}(1 - \beta)^2 ).
\eeq
\eel
\proof We only need to prove the lower bound.  From Lemmas \ref{lem:sad-1} and \ref{lem:msc-behaviour}we see that there is a $C_\eps >0$ so that for $N^{2/3}(E-2 ) \geq C_\eps$ we have with probability at least $1-\eps$ that
\beq
\frac{1}{N} \sum_{j=1}^N \frac{1}{ E - \lambda_j } \geq -\msc - N^{-1/3} \geq \beta + (1- \beta) - C\sqrt{(E - 2 ) } -N^{-1/3}.
\eeq
Hence if $N^{1/3}(1-\beta )\geq 10C (C_\eps)^2 + 10$ we see that $-m_N (E) >  \beta$ for $N^{2/3} (E - 2) \leq c_\eps (\beta-1)^2$.  This yields the claim, together with Proposition \ref{prop:sad-est}, which gives us the lower bound in the case that $N^{1/3} (1 - \beta) \leq10C (C_\eps)^2 + 10$.   \qed

Now we need the following.
\bel \label{lem:tilgam}
Let $\eps >0$ and $\xi >0$.  Assume $\beta \leq 1$. There is a constant $C_\eps$ on which the following holds with probability at least $1-\eps$ that, for large enough $N$,
\begin{align}
N^{2/3} | \tilgam - \gamma | \leq &C_\eps \bigg\{\frac{ \log^2 ( 1+  N^{1/3} ( 1 - \beta ) )}{ 1 + N^{1/3} (1 - \beta ) } \notag\\
&+ N^{-2/3+\xi} (1+ N^{1/3} (1- \beta) ) + \1_{ \{ N^{2/3} (1 - \beta)^2 \geq N^{1/15-\xi} \}} \frac{N^{\xi}}{1 + N^{1/3} ( \beta-1) } \bigg\}
\end{align}
\eel
\proof First we see that for $N^{2/3} ( \tilgam -2 ) \geq C_\eps$ we have with probability at least $1-\eps$ that,
\beq
N^{1/3}\left|m_N ( \tilgam) - \msc ( \tilgam) ) \right| \leq  C_\eps \left(  \frac{ \log^2 (N^{2/3} (\tilgam -2 ) )}{ N^{2/3} (\tilgam - 2)} +N^{-2/3+\xi} + \1_{ \{N^{2/3}( \tilgam -2) > N^{1/15-\xi} \}} \frac{N^{\xi}}{N^{2/3} (\tilgam -2 )} \right)
\eeq
for any $\xi >0$ and $N$ large enough, by Lemma \ref{lem:sad-1}. 
On the other hand, on $\GKA$ and if $N^{2/3}  ( \tilgam - 2) \geq 2 A$ we have,
\beq
| m_N ( \gamma ) - m_N ( \tilgam ) | \geq | \gamma - \tilgam | \frac{1}{N} \sum_{j=2}^N \frac{1}{ ( \lambda_j - \max\{ \tilgam, \gamma \} )^2}.
\eeq
There is an event of probability at least $1-\eps$ on which $\max\{ \tilgam, \gamma\} -2 \leq N^{-2/3} C_\eps ( 1 +  N^{2/3} ( 1 - \beta)^2 )$.  When this estimate holds and on $\GKA \cap \Fxi$ we have that for $j^{2/3} \leq N^{2/3} (1 _+N^{2/3} (1 - \beta)^2 )$ that 
\beq
N^{2/3} | \lambda_j - \max\{ \tilgam, \gamma \} | \leq K^{2/3} +  C_\eps (1 _+N^{2/3} (1 - \beta)^2 ).
\eeq
Hence, on this event,
\beq
\frac{1}{N^{4/3}} \sum_{j=2}^N \frac{1}{ ( \lambda_j - \max\{ \tilgam, \gamma \} )^2} \geq \frac{c_{K, \eps}}{ N^{1/3} (1 - \beta )+1 }.
\eeq
This yields the claim after taking $K, A$ large enough, from the equality
\beq
m_N ( \tilgam )- \msc( \tilgam) = m_N ( \tilgam) - m_N ( \gamma),
\eeq
which holds by definition of $\gamma, \tilgam$. \qed

Define now,
\beq
\tilgams = \max\{ \tilgam, 2 + sN^{-2/3} \}.
\eeq
We have the following. 
\bel
Let $\eps >0$.  Assume $\beta \leq 1$. There are $C_1$ and $C_2$ depending on $\eps$ so that any $s \geq C_1$ the following holds with probability at least $1-\eps$. 
\beq
| G( \gamma) - G ( \tilgams ) | \leq \frac{C_2(1+s^2)}{N}.
\eeq
\eel
\proof Assume $\GKA$ holds and choose $A$, $K$ large enough so that $\pp [ \GKA] > 1- \eps$.  Choose $C_1 > 2 A$.   We may further assume that the estimates of the previous two lemmas hold.

 By a second order Taylor expansion around $\gamma$ (using that $G' ( \gamma ) = 0$) we have,
\begin{align}
| G( \gamma ) - G ( \tilgams ) |  \leq C | \gamma - \tilgams |^2 \max_{x \in [ \gamma, \tilgams] } |G''(x) |.
\end{align}
On the event $\GKA \cap \Fxi$ for $\xi < 1/100$,  we have for $x > \lambda_1$,
\beq
\frac{1}{N^{1/3}} |G''(x) | \leq \frac{K}{N^{4/3} (x - \lambda_1 )^2} +C  \sum_{j > 10 K } \frac{1}{ j^{4/3} + N^{4/3}(x - \lambda_1)^2 } \leq \frac{K}{ N^{4/3} (x- \lambda_1)^2} + \frac{C}{ \sqrt{ N^{2/3} ( x- \lambda_1 )}}.
\eeq
This yields the claim after we note that we have,
\beq
| \gamma - \tilgams |^2 \leq C_\eps \frac{1+s^2}{N^{4/3}},
\eeq
on the event of Lemma \ref{lem:tilgam}.  \qed

\bep \label{prop:ht-G-asymp}
Let $\eps >0$.  Let $  \beta \leq 1$ and assume $(1-\beta)N^{1/3} \leq N^{1/10}$.  There is a $Q>0$ and $P>0$ so that there is an event with probability at least $1- \eps$ on which,
\beq
\left| G( \gamma ) - \left( R(2 + QN^{-2/3} ) + \beta^2/2 -1/2 + \log(\beta) \right) \right| \leq P \frac{1+ \log^5 (N^{1/3} (1-\beta) )}{N}
\eeq
where $R(z) := z  - N^{-1} \sum_i \log | z- \lambda_i |$.
\eep
\proof  First assume that the estimate of the previous lemma holds for some $\tilgams$.  We can assume also that $\GKA \cap \Fxi$ holds and $s > 2A$.   We take $Q= s$.  If $\tilgams = 2 + sN^{-2/3}$ we are finished, after noting that this implies that $1-\beta \leq CN^{-1/3}$, and that
\beq
\beta (2 + sN^{-2/3} ) - \beta^2 + 1/2 - \log ( \beta) = 2 + sN^{-2/3} +  \O (N^{-1})
\eeq 
in this case.  
 Otherwise, we integrate,
\begin{align}
\frac{1}{N} \sum_{j=1}^N \log ( \tilgam - \lambda_j) - \log (2 + sN^{-2/3} - \lambda_j ) &= -\int_{2 + sN^{-2/3}}^{\tilgam} m_N (t) \d t \notag\\
&= - \int_{2 + sN^{-2/3}}^{\tilgam} \msc (t) \d t -  \int_{2 + sN^{-2/3}}^{\tilgam} (m_N (t)- \msc (t)) \d t
\end{align}
By an application of Lemma \ref{lem:sad-1} we have,
\beq
\ee\left[ \1_{ \GKA \cap \Fxi } \left| \int_{2 + sN^{-2/3}}^{\tilgam} (m_N (t)- \msc (t)) \d t \right| \right] \leq C\frac{\log^3 (1+ N^{1/3} ( 1 - \beta ) )+N^{\xi} (1-\beta)^2 + }{N}.
\eeq
We can assume $\xi < 1/100$, so that the second term in the numerator is $\O(1)$.  
Hence by Markov's inequality this is less than $C_\eps \log^3 (1+ N^{1/3} ( 1- \beta ))$ with probability at least $1- \eps$.   Moreover, 
\beq
\int_{2 + sN^{-2/3}}^{ \tilgam } \msc (t) \d t = \int_2^{\tilgam} \msc(t) \d t +sN^{-2/3} +  \O (sN^{-1}).
\eeq
By explicit calculation,
\begin{align}
\int_{2}^{ \tilgam} \msc (t) \d t = \frac{1}{2} \left( 2-\tilgam^2/2 + \frac{ \tilgam}{2} \sqrt{ \tilgam^2-4} -2 \log ( ( \sqrt{\tilgam^2-4} + \tilgam)) + 2 \log(2) \right)
\end{align}
Since $\tilgam = ( \beta^2 +1 ) / \beta$ this simplifies to,
\begin{align}
\int_{2}^{ \tilgam} \msc (t) \d t  = \frac{1}{2} ( 1 - \beta^2 ) + \log ( \beta),
\end{align}
which yields the claim. 
\qed

\subsection{Low temperature case}

We now turn to the low temperature case.    The following is clear.
\bel \label{lem:sad-lt}
On $\JD$ we have for $E > \lambda_1$,
\beq 
\frac{1}{N} \sum_{i=1}^N \frac{1}{ E - \lambda_i } \leq 1 + D N^{-1/3} + \frac{1}{N^{1/3}} \frac{1}{ N^{2/3} (E - \lambda_1 )}.\,
\eeq
and so 
\beq
N^{2/3} ( \gamma - \lambda_1 ) \leq \frac{1}{ (N^{1/3} ( \beta - 1 )  - D)_+ }.
\eeq
\eel
In the low temperature regime we define $\hatgam$ by
\beq
\hatgam = \lambda_1 + \frac{1}{N ( \beta -1) }.
\eeq
We prove the following detailing the quality of approximation of $\gamma$ by $\hatgam$. 
\bel
Assume $\beta \geq 1$.  
On $\JD \cap \Sb \cap \GKA \cap \Fxi$ for $\xi = 1/100$ we have,
\begin{align}
\frac{ N^{2/3} | \gamma - \hatgam | }{ N^{4/3} ( \hatgam - \lambda_1 ) ( \gamma- \lambda_1 )}  \leq D + N^{2/3} ( \gamma - \lambda_1) C ( K b^{-1} + 1 ).
\end{align}
Therefore,
\begin{align}
N^{2/3} | \gamma - \hatgam | \leq \left( \frac{1}{ (N^{1/3} ( \beta -1 )  - D )_+} \right)^2\left( D + N^{2/3} ( \gamma - \lambda_1) C ( K b^{-1} + 1 ) \right).
\end{align}
\eel
\proof We have,
\begin{align}
\frac{1}{N ( \gamma - \lambda_1 ) } - \frac{1}{ N ( \hatgam - \lambda_1 ) } &=1 + \frac{1}{N} \sum_{j=2}^N \frac{1}{ \lambda_j - \lambda_1 } \notag\\
&+  \frac{\gamma - \lambda_1}{N} \sum_{j=2}^N \frac{ 1}{ (\lambda_j - \gamma)( \lambda_j -\lambda_1 ) }.
\end{align}
On $\Sb \cap \GKA \cap \Fxi$ we have,
\begin{align}
\frac{1}{N^{4/3}} \sum_{j=2}^N  \frac{ 1}{ (\lambda_j - \gamma)( \lambda_j -\lambda_1 ) } \leq C \left( \frac{K}{b}  +1 \right).
\end{align}
This yields the claim. \qed

\bep \label{prop:lt-G-asymp-1}
Suppose that $\JD \cap \Sb \cap \GKA \cap \Fxi$ hold and that $N^{1/3} (\beta-1) \geq 2 D$.    Let $Q > 2A$.  Then, 
\begin{align}
 G ( \gamma ) &= ( \beta -1 ) ( \lambda_1+N^{-1} ( \beta-1)^{-1} ) + \frac{1}{N} \log ( N^{1/3} ( \beta -1 ) )  \notag\\
&+ \left( 2 + N^{-2/3} Q - \frac{1}{N} \sum_{j=1}^N \log (2 + N^{-2/3}Q - \lambda_j ) \right) + \O ( N^{-1} ).
\end{align}
\eep
\proof  
Taylor expanding around $\gamma$ and using $G'( \gamma ) = 0$ we have,
\begin{align}
| G( \gamma) - G ( \hatgam)| \leq | \gamma - \hatgam|^2 \max_{ x \in [ \gamma, \hatgam]} |G''( x)|.
\end{align}
We see that on the events we are considering,
\beq
|G''(x) | \leq \frac{CK}{N} \left( \frac{1}{ (\lambda_1 - \hatgam)^2} + \frac{1}{ (\lambda_1 - \gamma)^2} \right) + CN^{1/3} \left( \frac{1}{ \sqrt{N^{2/3} ( \gamma-\lambda_1)}} + \frac{1}{ \sqrt{N^{2/3} ( \tilgam -\lambda_1 ) }} \right).
\eeq
Hence by applying the previous two lemmas and the assumption that $N^{1/3} ( \beta-1) \geq 2 D$ we have,
\begin{align}
|G( \gamma) - G ( \hatgam) | \leq \frac{C}{N( N^{1/3} ( \beta - 1) )^2}.
\end{align}
We now write
\beq
G ( \hatgam ) = ( \beta -1) \hatgam - \frac{1}{N} \log ( \hatgam - \lambda_1) + \left( \hatgam - \frac{1}{N} \sum_{j=2}^N \log ( \hatgam - \lambda_2 ) \right).
\eeq
Fix now a $Q > 2A$.  We have, using the second estimate of  Lemma \ref{lem:gpp} and that $|\tilgam - \lambda_1| \leq CN^{-2/3}$ and the assumption that we are working on the event $\JD$ that
\beq
\left| 1 + \frac{1}{N} \sum_{j=2}^N \frac{1}{\lambda_j - \hatgam } \right| \leq  CN^{-1/3},
\eeq
and so by a second order Taylor expansion around $\hatgam$ (applying the second estimate of Lemma \ref{lem:gpp}) we find,
\begin{align}
\left| \hatgam - \frac{1}{N} \sum_{j=2}^N \log ( \hatgam - \lambda_2 ) - (2+N^{-2/3}Q) + \frac{1}{N} \sum_{j=2}^N \log (  (2+N^{-2/3}Q)- \lambda_2 ) \right| \leq \frac{C}{N}(Q+A+1)^2
\end{align}
Also, 
\beq
\frac{1}{N} \log (2+QN^{-1/3} - \lambda_1 ) - \frac{1}{N} \log ( \hatgam - \lambda_1 ) = \frac{1}{N} \log ( N^{1/3} ( \beta-1) ) + N^{-1} \O ( \log (Q+A) ).
\eeq
This yields the claim. \qed

We also have the following expansion which will be used for the case that $N^{1/3} ( \beta-1)$ is bounded. 
\bep \label{prop:lt-G-asymp-2}
Assume that $\beta \geq 1$.  For every $\eps >0$ there is a $C_\eps >0$ so that the following holds with probability at least $1-\eps$ for $Q \geq C_\eps$.  We have for all $N$ large enough,
\beq
|G ( \gamma) - G(2 + Q N^{-2/3} ) | \leq \frac{C_\eps}{N} \left( 1 + N^{1/3} (\beta-1)\right)^2
\eeq
\eep
\proof This follows immediately from a second order Taylor expansion of $G$ around $\gamma$, the lower bound of Proposition \ref{prop:sad-est} and the first estimate of Lemma \ref{lem:gpp}. 
\qed

\section{Estimates on contour integral}
\label{sec:steep}

Recall that the saddle $\gamma$ is the unique solution $\gamma > \lambda_1$ to the equation
\beq
G' (\gamma ) = 0.
\eeq

We define the contour $\Gamma$ as the steepest descent curve, i.e., the curve in the complex plane passing through $\gamma$ and satisfying
\beq
\Im [ G ( z) ] = 0.
\eeq
Note that since $\bar{G(z) } = G ( \bar{z} )$ this curve is symmetric in the upper and lower half-planes and we will often just discuss its behavior in the upper half-plane.  Recall the following elementary consequences of the Cauchy-Riemann equations,
\beq
\partial_E \Re[f] = \Re[f'], \quad \partial_\eta \Re[ f'] = - \Im [ f'] , \quad \partial_E \Im [ f] = \Im [f'], \quad \partial_\eta \Im [f] = \Re [f'].
\eeq
for analytic $f$.  From the fact that $\partial_E \Im [ G (z) ] = \Im [ m_N ] > 0$ for $\eta >0$ we see that $\Im[G(z)]$ is a monotonic function of $E$, and so the curve maybe parameterized by $\eta$, and that once the curve leaves the real axis at $\gamma$ it never intersects the real axis again.  From the fact that $G'(z) \neq 0$ off of the real-axis we see that $\Gamma$ is $C^1$ curve by the implicit function theorem.  As $\eta \to \pi / \beta$, we see that $E \to - \infty$.  

We first prove,
\bel \label{lem:def-1} If $N \geq 20$ then for any $ a >\lambda_1$,
\beq
\int_{ a- \i \infty}^{a+ \i \infty} \e^{ N G(z)/2} \d z = \int_{\Gamma} \e^{ N G (z) /2 } \d z.
\eeq
\eel
\proof  Let $L = \sup_i | \lambda_i |$.   Let $C>0$. If $ \eta$ is sufficiently large we use that $\partial_\eta \Re[ G(z) ] \leq - ( 2 \eta)^{-1}$ for $|E| \leq 10L + 100(1+\beta) + C(1+|a|)$.   For $R$ sufficiently large depending on $L$ and $E \leq 0$ and $|z| \geq R$,
\beq
\Re[ G' (z) ] \leq \beta E - \log(|z|/2) \leq - \log(|z|/2).
\eeq 
Hence, if $N \geq 20$ we have that on the circle $|z| =R$ and $E \leq C$ for any $C>0$ we have for $R$ large enough,
\beq
| \e^{ N G(z)/2}|  \leq C' R^{-10}.
\eeq
Hence the contribution from this arc goes to $0$ as $N \to \infty$ and we can shift the contour to $\Gamma$. \qed

\bel\label{lem:cont-weak-upper}
Suppose that $\GKA \cap \Fxi$ hold with $\xi = 1/100$ and that $ ( \gamma - \lambda_1 ) \leq C' $.  Then, for $N$ sufficiently large,
\beq
\int_{-\infty}^\infty \left| \e^{ N (G (\gamma + \i t) - G ( \gamma ) )/2} \right| \d t \leq  C( N^{-2/3} (K+A+1) +  (\gamma - 2)_+)
\eeq
\eel
\proof We write,
\beq
\Im[ m_N (z ) ] \geq \frac{1}{N} \sum_{j \leq N \eta^{3/2} } \frac{ \eta}{ \eta^2 + ( \lambda_j - \gamma)^2 }.
\eeq
For $N$ sufficiently large we may assume that $K + N^{2/3} ( \gamma - 2)_+ + A \leq (2 C'+1) N^{2/3}$.  
If $(4 C'+2) N^{2/3} \geq N^{2/3} \eta \geq K + N^{2/3} ( \gamma - 2)_+ + A$, then on $\GKA \cap \Fxi$ the denominator in the sum above is smaller than
\beq
\eta^2 + ( \lambda_j - \gamma)^2 \leq C ( \eta^2  ).
\eeq
Hence, since $\Re[ G ( \gamma + \i t) ]$ is decreasing for $t>0$, we have
\beq
\Re[ G ( \gamma + \i t ) ] \leq - c \int_{N^{-2/3} (K+A) +  (\gamma - 2)_+}^t \eta^{1/2} \d \eta \leq - c t^{3/2}
\eeq
if $10 C' +10 > t > 2( N^{-2/3} (K+A) +  (\gamma - 2)_+)$.  Hence,
\beq
\int_{|t| \leq 1} \left| \e^{ N (G (\gamma + \i t) - G ( \gamma ) )/2} \right| \d t  \leq C( N^{-2/3} (K+A+1) +  (\gamma - 2)_+)
\eeq
For the regime $10(C'+1) \leq |t| \leq N^{10}$ we can just the fact that 
\beq
\Re[G(\gamma + \i t ) ] - \Re[ (G \gamma ) ] \leq \Re[G(\gamma + \i 10(C'+1) ) ] - \Re[ (G \gamma ) ] \leq -c,
\eeq
 by our above calculation, so this part of the integral is exponentially small.
Finally we use that on $\Fxi$ for $\eta \geq 1$ we have $\Im [ m_N ] \geq c / \eta$ and so,
\beq
\Re[ G' ( \gamma + \i t ) ] \leq - c \log ( t)
\eeq
for $|t| \geq 1$.  This allows us to estimate the rest of the integral.   \qed



We now move towards proving a lower bound.

\bel \label{lem:Gder-bds}
Suppose that the events $\GKA \cap \Fxi \cap \Sb$ hold with $\xi =1/100$ and $K, A, b>0$ and $K \geq A^{3/2}$ sufficiently large.  Let $\delta = N^{2/3} ( \gamma - \lambda_1)$.   Then for $N^{2/3} |z- \gamma| \leq \delta/2$ we have for $k \geq 2$, for a universal $C_1>0$,
\beq
\frac{N}{N^{2k/3}} |G^{(k)} (z) | \leq C_1^k (k-1)! \left[ \frac{1}{ \delta^k} + \frac{K}{ ( \delta + b)^k} + \frac{1}{ \delta^{k-3/2}} \right],
\eeq
Assume further that $\delta \leq C' N^{2/3}$ for some $C'>0$.   Then for $k=2, 3, 4$,
\beq
\frac{N}{N^{2k/3}} |G^{(k)} (\gamma) | \geq c \left( \frac{1}{ \delta^k } + \1_{ \{ \delta^{3/2} > K \} } \frac{1}{ \delta^{k-3/2}} \right),
\eeq
where $c$ depends on $C'$. 
\eel
\proof For the upper bound we have,
\beq
\frac{N}{N^{2/3k}} |G^{(k)} (z) | \leq C^k (k-1) \left[ \frac{1}{ \delta^k} + \frac{K}{ ( \delta+b)^k} + \sum_{j > K } \frac{1}{ (N^{2/3} ( \lambda_j -z))^k} \right].
\eeq
Under our assumptions we have that for $K$ large enough and $j > K$ on $\GKA \cap \Fxi$,
\beq
N^{2/3} |\lambda_j - z | \leq C ( \delta  + j^{2/3} ).
\eeq
Hence the estimate follows from breaking up the sum into $j^{2/3} \leq \delta$ and $j^{2/3} \geq \delta$.  

For the lower bound,
\beq
\frac{N}{N^{2k/3}} |G^{(k)} (\gamma) | \geq  \frac{c}{ \delta^k} + \sum_{j^{2/3} \leq \delta} \frac{1}{ ( N^{2/3} ( \lambda_j - \gamma ))^k }.
\eeq
As long as $\delta \geq K^{2/3}$ the factor in the denominator in the second term is larger than $c \delta$.  Our assumed upper bound on $\delta$ ensures that the sum has at least $(\delta/C')^{3/2}$ terms.  This yields the lower bound. \qed

With this, we now prove the following concerning the behavior of the contour $\Gamma$ in a small disc around the saddle $\gamma$.

\bep \label{prop:Gam-est}
Suppose that the events $\GKA \cap \Fxi$ hold with $\xi = 1/100$.  Define $\delta = N^{2/3} ( \gamma - \lambda_1)$ and assume that $\delta \leq C_2 N^{2/3}$.  There are constants $c_1, c_2$ and $C_1$ depending on $K$ and $C_2$ so that for $N^{2/3} |z - \gamma | \leq c_1 \delta$ we have for $z \in \Gamma$,
\beq
c_2 N^{2/3} |E-\gamma| \delta \leq N^{4/3} \eta^2 \leq C_1 N^{2/3} |E-\gamma| \delta .
\eeq
\eep
\proof Throughout the proof we will always assume $\GKA \cap \Fxi$ holds without comment.  
Let us change variables to,
\beq
N^{2/3} (z - \gamma ) = w=  u + \i v.
\eeq
By either a power series expansion (using estimates from the previous lemma) or Taylor's theorem with explicit remainder we have for $|w| \leq c\delta $ that,
\beq \label{eqn:c-1}
N (G (\gamma + z ) - G ( \gamma )) = w^2 T_2 + w^3 T_3 + w^4 g (w)
\eeq
where,
\begin{align}
c' \left( \frac{1}{ \delta^k} + \1_{ \{ \delta^{3/2} > K \} } \frac{1}{ \delta^{k-3/2} } \right) \leq  |T_k|  \leq C \left[ \frac{K}{ \delta^k} + \frac{1}{ \delta^{k-3/2}} \right],
\end{align}
and,
\beq
|g(w) | \leq  C \left[ \frac{K}{ \delta^4} + \frac{1}{ \delta^{5/2}} \right].
\eeq
Furthermore, $T_2 >0$ and $T_3 < 0$.  Along the real axis $g(w)$ and its derivatives are real so,
\beq
\left| \Im [ g(w) ] \right| \leq C |v| \left( \delta^{-5} + \delta^{-7/2} \right).
\eeq
Taking the imaginary part of \eqref{eqn:c-1}  and setting it equal to $0$ we obtain
\beq
0 =  2v u T_2 +  u^2 v 3 T_3 - v^3 T_3 + (4 u^3 v - 4 u v^3 ) \Re[g (w) ] + (u^4 - 6 v^2 u^2 + v^4 ) \Im [ g(w) ].
\eeq
Dividing through by $v$ we find,
\beq
v^2 T_3 = 2 u T_2 + u^2 3 T_3 + (4 u^3 - 4 u v^2) \Re[g (w) ] + (u^4 - 6 v^2 u^2 + v^4 ) (\Im [ g(w)]/v).
\eeq
Note that we see that
\beq
\frac{K + \delta^{3/2}}{\delta^k} \frac{c}{ K+1} \leq |T_k| \leq C \frac{ K+ \delta^{3/2}}{\delta^k}.
\eeq
Hence, for constants $c_1, c_2$ and $C_1$ depending on $K$ we have that for $|w| \leq c_1 \delta$ that
\beq
c_2 |u| \delta \leq v^2 \leq C_1 |u| \delta.
\eeq
This yields the claim after noting the signs of $T_2$ and $T_3$. \qed
\bep \label{prop:cont-lower}
Suppose that $\GKA \cap \Fxi$ with $\xi =1/100$ and that $\delta  = N^{2/3} ( \gamma - \lambda_1)$.  Then there are constants $c, C$ depending on $K$ so that,
\beq
-\i \int_{ \Gamma} \e^{ N(G(z) - G( \gamma ) )/2 } \d z \geq c \delta N^{-2/3} \e^{ - C \sqrt{ \delta}}.
\eeq
\eep
\proof The integral is real so for any parameterization of $\Gamma$ by $z = \d x + \i \d y$, only the $\d y$ part contributes. Moreover, the integrand is positive and along $\Gamma$, the imaginary part of $ z \in \Gamma$ is monotonically increasing.  Hence  for a lower bound we can take the portion of the curve intersecting a disc of radius $c_1 N^{-2/3} \delta$ centered at $\gamma$, where $c_1$ is from Proposition \ref{prop:Gam-est}.  From Proposition \ref{prop:Gam-est}, the curve exits the disc at a point with $\Im[z] \geq c \delta N^{-2/3}$.   From the fact that $G'( \gamma ) = 0$ and the upper bounds on $G''(z)$ of Lemma \ref{lem:Gder-bds} we get
\beq
|\Re[ G(z) ] | \leq C \delta^{1/2}
\eeq
for $z$ in this disc.  This yields the claim. \qed

\subsection{Better upper bound for contour integral in low temperature regime}

In the low temperature regime we can obtain a better upper bound.  Recall the definition of the contour $\Gamma$.  For $a >0$ we define,
\beq
\Gama := \{ z \in \Gamma : E - \lambda_1 \geq -aN^{-2/3} \}.
\eeq
and we also introduce the shorthand,
\beq
\mfb = \frac{1}{N^{1/3} ( \beta - 1 ) }.
\eeq
The next proposition and lemma control the behavior of the contour as well as the integrand along the contour.
\bep \label{prop:lt-int-1}
Let $\eps >0$.  There is an $a >0$ depending on $\eps$ and a $C_\eps >0$ so that the following holds with probability at least $1-\eps$.  If $N^{1/3} ( \beta -1 ) \geq C_\eps$ and $ c_1>0$ is a constant we have for $z \in \Gama$ satisfying $E - \gamma \leq -c _1\mfb N^{-2/3}$ that there is a $C >0$ depending on $c_1$ so that,
\beq
\frac{1}{C C_\eps } \mfb \leq N^{2/3} \eta \leq \mfb C C_\eps,
\eeq
with probability at least $1-\eps$.   
Moreover, for $z \in \Gama$, $\eta$ is a monotonic function of $E$ and
\beq
\Re[ G'(z) ] \geq \frac{1}{2} ( \beta-1) \1_{ \{ E \leq \lambda_1 \}}.
\eeq
\eep
\proof Since $\partial_E \Im [ G(z) ] = \Im [ G'(z) ] >0$ we see that the contour intersects each horizontal line the complex plane only once.  Hence, the lower bound on $\eta$ follows from Propositions \ref{prop:sad-est} and  \ref{prop:Gam-est}.   

For the upper bound we will argue by finding a vertical line $\eta \to E_* + \i \eta$ on which the sign of $\Im[G(z)]$ changes at a point $\eta N^{2/3} \leq C \mfb$.

Let us assume that $\Sb$ holds, taking $b$ small enough so that $\pp [ \Sb] \geq 1- \eps$.  By taking $\mfb$ small we can assume that 
\beq
N^{2/3} ( \gamma - \lambda_1 ) \leq C_\eps \mfb \leq b /100
\eeq
 with probability at least $1-\eps$ by applying the upper bound of Lemma \ref{lem:sad-lt}. 
For $|z- \lambda_1 | \leq b /10$ we have the following expansion,
\begin{align} \label{eqn:c-2}
G(z) = ( \beta + \Xi ) ( z - \lambda_1 ) + \beta \lambda_1 - \frac{1}{N} \log ( z - \lambda_1 ) + (z-\lambda_1)^2 g(z)
\end{align}
for a function $g(z)$ satisfying $|g(z) | \leq N^{1/3} C$ on $\GKA \cap \Fxi$ (and $\xi = 1/100$) where $C$ depends on $K,  \eps, b$ and  we defined $\Xi$ as
\beq
\Xi := \frac{1}{N} \sum_{j > 1}\frac{1}{ \lambda_j - \lambda_1 } .
\eeq
For $E < \lambda_1$ we have $\Im [ G(E + \i 0) ] <0$ (as a boundary value).   On $\JD$  and the current event we are considering, we can derive from \eqref{eqn:c-2},
\beq
 \Im [ G(z ) ] \geq (\beta-1 - D N^{-1/3} ) \eta - \frac{\pi}{N} - CN^{1/3} ( (E- \lambda_1 )^2+ \eta^2 ).
\eeq
We choose $E_* < \lambda_1$ so that  $|E_* - \lambda_1 | = \min\{ c N^{-2/3}, bN^{-2/3}/100\}$.  Assume that $(\beta-1) \geq 2 D N^{-1/3}$  so that with $y = N (\beta-1) \eta$ we have,
\beq
N \Im[ G(E_* + \i \eta) ] \geq \frac{y}{2} - 10 - C y^2  \mfb^2.
\eeq
If $\mfb$ is small enough this will be positive for $y = 100$, or for $\eta = 100N^{-2/3}  \mfb$.  The fact that we used a Taylor expansion above means that we need to ensure that this choice of $y$ does not conflict with the constraint that $N^{2/3} \eta \leq b /100$.  So we choose $\mfb$ small enough so that this is the case.  This means that for this chosen $E_* < \lambda_1$ the equation $\Im [ G(E+ \i \eta) ] = 0$ has a solution that is less than $C N^{-2/3} \mfb$.  This yields the upper bound we desire, with the $a$ for $\Gama$ coming from $E_*$. 

For monotonicity, by the equation
\beq
\frac{ \d \eta}{\d E} = - \frac{ \Im [G' (z) ]}{ \Re[ G' (z) ]}
\eeq
it suffices to check that $\Re[G'(z) ] > 0$.  This argument is slightly different depending on whether $E$ is less than or greater than $\lambda_1$.  For $\lambda_1 < E < \gamma$ consider first as a function of $\eta$,
\beq
\Im [G(z) ] = \beta \eta - \frac{1}{N} \sum_{j} \log (E + \i \eta - \lambda_j ),
\eeq
and
\beq
\partial_\eta \Im [ G(z) ] = \beta + \frac{1}{N} \sum_{j=1}^N \frac{ \lambda_j - E }{ ( \lambda_j - E)^2 + \eta^2 }.
\eeq
This is a convex function of $\eta$ that is $0$ at $\eta = 0$ and $\del_\eta \Im [G (z) ] \vert_{\eta = 0} < 0$ by definition of the saddle $\gamma$.  We also see that for large $\eta$ that $\Im [ G(z) ] >0$.  The positive solution to $\Im [ G(z) ] = 0$ is the point on the contour $\Gamma$ lying above $E$ and we see, by the strict convexity, that $\del_\eta \Im [G (z) ] = \Re [ G' (z) ] > 0$ at this point.  This completes the case $\lambda_1 < E < \gamma$.  

For $E < \lambda_1$ and $z \in \Gama$ we have, on $\JD \cap \Sb \cap \GKA \cap \Fxi$, (choose all these parameters so that the events  hold with probability at least $1-\eps$)
\begin{align}
\Re[ G'(z) ] \geq \beta -1 + \frac{1}{N} \frac{ \lambda_1 - E}{(\lambda_1 - E)^2 + \eta^2 } - DN^{-1/3} - |z-\lambda_1 |  CN^{1/3}
\end{align}
using a similar expansion as before, as long as $|z- \lambda_1| \leq b/100$.  Here, $C$ depends on $\eps$ but not on $\mfb$.  We can decrease $a$ and assume that $\mfb$ is sufficiently small so that the portion of $\Gama$ for which $E \leq \lambda_1$ lies in the disc $|z-\lambda_1| \leq b /100$ with probability at least $1-\eps$, using the upper bound just proven for $\eta$.  We can further let $a$ and $\mfb$ be small enough so that the very last term in the above inequality is less than $|z-\lambda_1|C N^{1/3} \leq  N^{-1/3}$.    Now, we see that if $(\beta-1)N^{1/3}$ is large enough, depending on $D$ this is positive and moreover exceeds $(\beta-1)/2$. \qed

\bel
Under the conditions of the previous Proposition, we have that for $z \in \Gama$ that,
\beq
\Re[G(z) - G ( \gamma ) ] \leq -c(\beta-1) (E- \lambda_1 ) \1_{  \{ E \leq \lambda_1 \} }
\eeq
and if $E \leq \lambda_1$ that,
\beq
\left| \frac{ \d \eta}{ \d E} \right| \leq C
\eeq
\eel
\proof We have,
\beq
\frac{\d}{\d E} \Re[ G (z) ] = \Re[ G'(z) ] + \frac{ ( \Im [ G'(z) ] )^2}{ \Re[ G'(z) ]} \geq \Re[ G' (z) ].
\eeq
The first claim is then an application of the lower bound for $\Re[G'(z) ]$ we found in the previous proposition.  For the second we see that for $E \leq \lambda_1$ 
\beq
\Im [G'(z) ] \leq \frac{1}{N \eta} + C \eta N^{1/3} \leq  C ( \beta-1),
\eeq
on $\GKA \cap \Fxi$. This together with the equality
\beq
\left| \frac{ \d \eta}{ \d E} \right| = \left| \frac { \Im [G' (z) ]}{ \Re[ G' (z) ] } \right|
\eeq 
and the lower bound for $\Re[ G'(z)]$ yields the claim.
\qed

With the two previous results we can now estimate the integral.
\bep \label{prop:lt-int-upper}
Let $\eps >0$.  There is a $C_\eps$ so that if $N^{1/3} ( \beta-1) \geq C_\eps$ then with probability at least $1-\eps$,
\beq
\left| \frac{1}{ 2 \pi \i } \int_\Gamma \e^{ N (G(z) - G( \gamma ) )/2 } \d z \right| \leq C_\eps \mfb N^{-2/3}.
\eeq
\eep
\proof For any $a$ we can deform the contour to be the union of $\Gama$ and the vertical lines extending from the end of $\Gama$ to $ \pm \i \infty$, using the argument of Lemma \ref{lem:def-1}.  First we estimate the contribution of $\Gama$.  First, for $ \lambda_1 \leq E \leq \gamma$ we just use that the length of this contour is $\O ( \mfb N^{-2/3})$ and that $\Re[ G(z)]$ is decreasing there, so that the integrand is bounded by $1$.  The estimate on the length of the contour follows from the fact that $\eta$ and $E$ are monotonically related, the upper bound on $\eta$ from Proposition \ref{prop:lt-int-1}, and the upper bound on the saddle of Lemma \ref{lem:sad-lt}.   

For the rest of $\Gama$ we can assume that $(\d \eta ) / ( \d E)$ is bounded due to the previous lemma and so the contribution here is bounded by,
\beq
\int_{0}^\infty \exp \left[ - c N ( \beta-1) t \right] \d t \leq C \frac{ \mfb}{N^{2/3}},
\eeq
where we also used the previous lemma to estimate the integrand. 

The remaining integral consists of the vertical lines from the endpoints of $\Gama$ to  $\pm \i \infty$.  Their contribution is bounded by,
\beq
C \e^{ -c ( \beta-1)N^{1/3} }\left(  \int_{ \eta >   N^{-2/3} } \exp\left[ - N \int_{N^{-2/3}}^\eta \Im [G' (E_a + \i t ) ] \d t \right] \d \eta +N^{-2/3} \right)
\eeq
where $E_a = \lambda_1 - a N^{-2/3}$.  Note that we estimated the contributino for $\mfb c N^{-2/3} \leq \eta \leq N^{-2/3}$ by $N^{-2/3}$.  The same argument from the start of the proof of  Lemma \ref{lem:cont-weak-upper} shows that $\Im [ G' (z) ] \geq c \eta^{1/2}$ for $N^{-2/3} \leq \eta \leq N^{-2/3+1/10}$.   Hence, the contribution in this range is $\O (N^{-2/3} )$.  At the point $\eta = N^{-2/3+1/10}$ the integrand is exponentially small, so we can do away with the contribution for $N^{-2/3+1/10} \leq \eta \leq N^C$.  For $\eta \geq 1$ we have $\Im [G'(z) ] \geq c / \eta$ so $\Re[ G(z) ] \leq -c \log ( \eta)$ here. Hence, the contribution for $ \eta > N^{C}$ is $\O(N^{-10})$ if $C$ is large enough.  This yields the claim. \qed

\subsection{Tighter estimates in high temperature regime}

We first have the following.
\bel \label{lem:ht-b1}
Suppose that $\GKA \cap \Fxi$ holds with $\xi = 1/100$.   Assume that there is a constant $C_1>0$ so that $\gamma -2 \leq C_1$ and furthermore that $N^{2/3} (\gamma - 2) \geq 10 (A+K)$.   There is a $C >0$ so that,
\beq
\frac{1}{ C( \gamma-2)^{1/2} } \leq |G'' ( \gamma ) | \leq \frac{C}{ ( \gamma-2)^{1/2}}, \qquad |G'''( \gamma ) | \leq \frac{C}{ ( \gamma-2)^{3/2}}
\eeq
and for $0 \leq \eta  \leq 1$,
\beq
\Im [ m_N ( \gamma + \i \eta ) ] \geq \frac{1}{C} \frac{ \eta}{ \sqrt{ \eta +  (\gamma-2) }}.
\eeq
\eel
\proof Under the assumptions we see that,
\beq
\frac{1}{C} ( N^{2/3}( \gamma - 2) + j^{2/3} ) \leq N^{2/3} (\gamma  - \lambda_j )\leq C( N^{2/3} ( \gamma - 2 ) + j^{2/3} ),
\eeq
and so 
\beq
\frac{1}{C} \sum_{j=1}^N \frac{1}{ (N^{2/3} ( \gamma - 2 ) + j^{2/3} )^k} \leq \sum_{j=1}^N \frac{1}{ ( N^{2/3} (\gamma  - \lambda_j ) )^k } \leq C \sum_{j=1}^N \frac{1}{ (N^{2/3} ( \gamma - 2 ) + j^{2/3} )^k }
\eeq
for $k=2, 3$.  The upper bound follows from consider separately the contributions from $j^{2/3} \leq N^{2/3}  ( \gamma-2)$ and $j^{2/3} \geq N^{2/3} ( \gamma-2)$ separately, and estimating the factor in the denominator from below by $N^{2/3} ( \gamma-2)$ and $ j^{2/3}$ in each case, respectively.  The lower bound comes from the terms $j^{2/3} \leq N^{2/3}  ( \gamma-2)$ for which the factor in the denominator is less than $2 N^{2/3} ( \gamma -2 )$.

For the lower bound on $\Im [ m_N ( \gamma  + \i \eta )$ we have, 
\beq
\Im [m_N ( E + \i \eta ) ] \geq c \frac{1}{N} \sum_{j=1}^N \frac{\eta}{ ( \gamma - 2)^2 + N^{-4/3} j^{4/3} +\eta^2} \geq \frac{c}{N} \sum_{ j^{2/3} \leq N^{2/3} ( \gamma -2 )+N^{2/3} \eta} \frac{ \eta}{ ( \gamma -2 )^2 + \eta^2}
\eeq
There are at least $c ( N^{2/3} ( \gamma - 2) + N^{2/3} ( \eta ) )^{3/2}$ terms in the final sum  on the RHS due to our assumption $N^{2/3} ( \gamma -2 ) \leq C_1$. This yields the claim. \qed

\bep \label{prop:ht-b1} Let $\eps >0$ and $ \beta \leq 1$. There is a $C_\eps >0$ so that if $N^{1/3} ( 1 - \beta ) \geq C_\eps$ then the following holds with probability at least $ 1- \eps$.
\beq
\frac{1}{C_\eps} \frac{ ( N^{1/3} (1 - \beta) )^{1/2} }{N^{2/3}} \leq \int_{\rr} \exp\left[ \frac{N}{2} ( G ( \gamma + \i t ) - G ( \gamma ) ) \right] \d t \leq C_\eps \frac{ ( N^{1/3} (1 - \beta) )^{1/2} }{N^{2/3}}
\eeq
\eep
\proof
We assume that $\GKA \cap \Fxi$ hold with probability at least $1- \eps$.  We furthermore assume that the event of Lemma \ref{lem:ht-saddle-bd} holds.  We then assume that $(1-\beta)N^{1/3}$ is sufficiently large so that the estimates of Lemma \ref{lem:ht-saddle-bd} imply that $\gamma - 2$ is sufficiently large that the previous lemma applies.  Fixing a small $c_1 >0$ we split the integral into a few regions
\begin{enumerate}[(i)]
\item $ |t| \leq 2 c_1 ( \gamma -2 )^{1/2} N^{-1/3}$ \label{it:1}
\item $2c_1 ( \gamma -2 )^{1/2} N^{-1/3}( \gamma -2 ) < |t| < ( \gamma - 2)$ \label{it:2}
\item $( \gamma -2 ) < |t| <1 $ \label{it:3}
\item $ |t| > 1 $ \label{it:4}
\end{enumerate}
For region \eqref{it:1} we use, the second order Taylor expansion,
\beq
\frac{N}{2} \left( G ( \gamma + \i t ) - G( \gamma ) \right ) = - \frac{ Nt^2}{4} G'' ( \gamma ) + \O\left( |t|^3 N ( \gamma - 2)^{-3/2} \right)
\eeq
where we used the estimates of Lemma \ref{lem:ht-b1} and the fact that $|G''' ( \gamma + \i t ) | \leq | G''' ( \gamma )|$.  Hence, we can choose $c_1$ sufficiently small so that,
\begin{align}
\exp\left[ \frac{N}{2} (G ( \gamma + \i t ) - G ( \gamma ) ) \right] = \exp\left[ - \frac{Nt^2}{4} G'' ( \gamma ) \right] \left(1 + f(t) \right)
\end{align}
for a function that satisfies $|f(t)| \leq \frac{1}{100}$ for $|t| \leq 100 c_1 (\gamma-2)^{1/3} N^{-1/3}$.  We see also that if $N^{1/3} ( \beta -1 )$ is large enough, then
\beq
\Re[ G ( \gamma + \i c_1 ( \gamma -2 )^{1/2} N^{-1/3} ) - G ( \gamma ) ] \leq -c N^{-2/3} ( \gamma -2)^{1/2}.
\eeq
 Hence,
\begin{align}
&\int_{ |t| \leq c_1 ( \gamma-2)^{1/3} N^{2/3} } \exp\left[ \frac{N}{2} (G( \gamma + \i t ) - G ( \gamma ) ) \right] \d t  \notag\\
=& \frac{( N^{2/3} ( \gamma - 2 ))^{1/4}}{ N^{2/3}} \frac{1}{ \sqrt{ 2 \pi G'' ( \gamma) ( \gamma -2)^{1/2} } } \left( 1 +\O \left( \e^{  -cN^{1/3} ( \gamma -2 )}  \right)\right) \left(1 + X \right)
\end{align}
where $|X| \leq \frac{1}{100}$.   For region \eqref{it:2}  we have, using that $\Re[ G ( \gamma + \i t)]$ is decreasing for $t>0$ that,
\begin{align}
\Re[ G( \gamma + \i t) -G ( \gamma ) ] \leq \Re[ G( \gamma + \i t) -G ( \gamma + \i c_1 ( \gamma -2 )^{1/2} N^{-1/3} ) ] - c N^{-2/3} ( \gamma -2)^{1/2}
\end{align}
and then,
\begin{align}
\Re[ G( \gamma + \i t) -G ( \gamma + \i c_1 ( \gamma -2 )^{1/2} N^{-1/3} ) ]  &= - \int_{ c_1 (\gamma-2)^{1/3} N^{-1/3}}^{t} \Im [m_N (\gamma + \i \eta ) ] \d \eta  \notag\\
&\leq -c  \int_{ c_1 (\gamma-2)^{1/3} N^{-1/3}}^{t} \eta / \sqrt{ ( \gamma -2 ) } \d \eta \notag\\
&\leq - c \frac{ t^2}{ \sqrt{ ( \gamma -2 ) }}
\end{align}
where we used Lemma \ref{lem:ht-b1} in the second line and that $t > 2 (c_1 (\gamma-2)^{1/3} N^{-1/3})$ in the third line.  Hence,
\begin{align}
& \int_{ 2c_1 ( \gamma -2 )^{1/2} N^{-1/3}( \gamma -2 ) < |t| < ( \gamma - 2)} \left| \exp\left[ - \frac{N}{2} (G( \gamma + \i t ) - G( \gamma ) ) \right] \right|\d t \notag\\
 \leq & C \e^{ - c N^{1/3} ( \gamma -2)^{1/3} } \int_{ |t| < (\gamma - 2) } \exp\left[ - c N t^2 (\gamma-2)^{-1/2} \right] \d t  \leq  C \e^{ - c N^{1/3} ( \gamma -2)^{1/3} }  \frac{ (N^{2/3} ( \gamma -2))^{1/4}}{N^{2/3}}.
\end{align}
For region \eqref{it:3} we use that 
\begin{align}
\Re[ G ( \gamma + \i t ) -G( \gamma ) ] =& \Re[ G ( \gamma + \i t ) - G( \gamma + i ( \gamma -2)/2 ) ] + \Re[ G( \gamma + \i (\gamma-2)/2 ) - \Re[ G ( \gamma ) ] \notag\\
 \leq& \Re[ G ( \gamma + \i t ) - G( \gamma + i ( \gamma -2)/2 ) ] - c ( \gamma -2)^{3/2},
\end{align}
and then,
\begin{align}
Re[ G ( \gamma + \i t ) - G( \gamma + i ( \gamma -2)/2 ) ] &= - \int_{ ( \gamma -2)/2}^{ t} \Im [m_N ( \gamma + \i \eta ) ] \d \eta \notag\\
&\leq - c \int_{ ( \gamma -2)/2}^{ t} \sqrt { \eta } \d \eta \leq - c t^{3/2}
\end{align}
where we used Lemma \ref{lem:ht-b1} and the fact that $t > (\gamma -2)$.  Hence,
\begin{align}
&\int_{ ( \gamma-2) < |t| < 1} \left| \exp\left[ \frac{N}{2} (G ( \gamma + \i t ) - G ( \gamma ) ) \right] \right| \d t \notag\\
\leq & C \e^{ - c N ( \gamma -2 )^{3/2} } \int_{|t| \leq 1} \exp \left[ - c N t^{3/2} \right] \d t \leq C N^{-2/3} \e^{ - c N ( \gamma -2 )^{3/2} }.
\end{align}
For $|t| > 1 $ we use,
\begin{align}
\Re[ G ( \gamma + \i t ) - G ( \gamma ) ] \leq \Re[ G ( \gamma + \i t ) - G ( \gamma + \i ) ]  - c,
\end{align}
and that for $ t > 1$ we have $\Im [ m_N ( \gamma + \i t ) ] \geq  c t^{-1}$ and so ,
\beq
\Re[ G ( \gamma + \i t ) - G ( \gamma ) ] \leq - c \log ( t) - c.
\eeq
Hence the contribution from region \eqref{it:4} is $\O (\e^{-N})$.   The leading order contribution from region \eqref{it:1} is positive and is,
\beq
 \frac{( N^{2/3} ( \gamma - 2 ))^{1/4}}{ N^{2/3}} \frac{1}{ \sqrt{ 2 \pi G'' ( \gamma) ( \gamma -2)^{1/2} } }.
\eeq
By taking $N^{1/3} ( \beta -1 )$ large enough we see that this ensures $(\gamma -2 )N^{2/3}$ large enough to ensure that all the other contributions from the integral are less than half  of this quantity.  This yields the claim. \qed

\section{Expansion of free energy}
\label{sec:expand}
Define
\beq
X_Q = \sum_j \log |2 + QN^{-2/3}  - \lambda_j |  - QN^{1/3} -N/2 + \log(N)/6.
\eeq

\subsection{High temperature regime}
\bep  \label{prop:fe-ht}
Assume $\beta \leq 1$.  
For any $\eps >0$ there is a $C_\eps >0$ so that with probability at least $1-\eps$ we have for a large $Q$ depending on $\eps >0$,
\begin{align}
\left| N(F_N - \beta^2/4) + \log(N)/12 +X_Q/2 \right| \leq C_\eps \left( 1 + N^{1/3} (\beta -1 ) \right).
\end{align}
\eep
\proof We first have by Lemma \ref{lem:rep} and \eqref{eqn:stir} that,
\begin{align}
N( F_N- \beta^2/4)  &= \frac{N}{2} G ( \gamma ) + \frac{1}{2} \log(N) - \frac{N}{2} \log ( \beta) - \frac{N}{2} - \beta^2 N/4 + \O(1) \notag\\
&+ \log \left( \frac{1}{2 \pi \i } \int \exp \left[ \frac{N}{2} (G(z) - G( \gamma ) )\right] \d z \right)
\end{align}
By Lemma \ref{lem:cont-weak-upper}, Proposition \ref{prop:cont-lower} and Lemma \ref{lem:ht-saddle-bd} we 
 see that there is a $C_\eps >0$ on which,
\beq \label{eqn:ht-b1} 
\left| \log \left( \frac{1}{2 \pi \i } \int \exp \left[ \frac{N}{2} (G(z) - G( \gamma ) )\right] \d z \right) + \frac{2}{3} \log(N) \right| \leq C_{\eps} \left( 1 + N^{1/3} ( 1-\beta) + \log (1 + N^{1/3} ( \beta -1 ) ) \right).
\eeq
with probability at least $1 - \eps$.  We have,
\begin{align}
 &\frac{N}{2} G ( \gamma ) + \frac{1}{2} \log(N) - \frac{N}{2} \log ( \beta) - \frac{N}{2} - \beta^2 N/4 \notag\\
 =& \frac{N}{2} \left( G ( \gamma ) - \beta^2/2 - \log (\beta) -1 \right) + \frac{1}{2} \log(N) \notag\\
=& \frac{N}{2} \left( G ( \gamma) - \beta^2/2 - \log ( \beta) - 1 + X_Q / N - \frac{1}{6} \log(N) / N \right) + \frac{7}{12} \log(N) - X_Q/2.
\end{align}
Now note that,
\beq 
G ( \gamma) - \beta^2/2 - \log ( \beta) - 1 + X_Q / N - \frac{1}{6} \log(N) / N = G( \gamma) -\beta^2/2 -R (2  +QN^{-2/3} ) + 1/2- \log ( \beta) ,
\eeq
where $R(z)$ is as in Proposition \ref{prop:ht-G-asymp}.  This proposition also gives us a bound for the RHS, and gives us the choice of $Q$.  This concludes the proof. \qed

In the above proof if we instead use Proposition \ref{prop:ht-b1} to estimate \eqref{eqn:ht-b1} we easily deduce the following.
\bep
\label{prop:ht-b2}  Let $\eps >0$ and $\beta \leq 1$.  There is a $C_\eps >0$ and a $Q>0$ so that if $N^{1/3} ( 1- \beta ) \geq C_\eps$ then
\beq
\left| N(F_N - \beta^2/4) + \log(N)/12 +X_Q/2  \right| \leq C_\eps (1 + \log^t (N^{1/3} ( 1 - \beta ) ) )
\eeq
\eep

\subsection{Low temperature regime}

\bep\label{prop:fe-lt}
For any $\eps >0$ there is a $Q>0$ and $C_\eps >0$ so that with probability with at least $1- \eps $ we have that,
\begin{align}
\bigg| N (F_N - [ \beta - &3/4 - \log(\beta)/2 ] )  \notag\\
- & \left( - X_Q +N (\beta-1) ( \lambda_1 - 2) - \log(N)/12 - \log(N^{1/3} ( \beta -1 ) + 1 ) /2\right) \bigg| \leq C_\eps
\end{align}
\eep
\proof  We first have, by combining Propositions \ref{prop:lt-G-asymp-1} and \ref{prop:lt-G-asymp-2} that, for any $\eps >0 $ there is a $Q>0$ and $C_\eps >0$ so that,
\beq
\left|G ( \gamma ) -\left( 2 \beta- \frac{1}{2} + ( \beta - 1) ( \lambda_1 - 2) + \frac{1}{N} \log ( N^{1/3} ( \beta -1) + 1 ) + \frac{ \log(N)}{6 N } - X_Q /N\right) \right| \leq \frac{C_\eps}{N},
\eeq
with probability at least $1-\eps$. 
Furthermore, by Propositions \ref{prop:cont-lower}, \ref{prop:sad-est} and \ref{prop:lt-int-upper} we see that there is a $C_\eps >0$ so that with probability at least $1-\eps$,
\beq
\left| \log \left( \frac{1}{ 2 \pi \i }  \int_{\Gamma} \exp \left[ N (G(z) - G(\gamma) )/2 \right] \d z \right) + \log (N^{2/3} ( N^{1/3} ( \beta -1 ) + 1 ) )  \right| \leq C_\eps.
\eeq
Hence,
\begin{align}
N F_N &= \frac{N}{2} G( \gamma) + \log \left( \frac{1}{ 2 \pi \i }  \int_{\Gamma} \exp \left[ N (G(z) - G(\gamma) )/2 \right] \d z \right) + \frac{1}{2} \log(N) - \frac{N}{2} \log ( \beta) - \frac{N}{2} + \O(1) \notag\\
&= N \left( \beta - \frac{3}{4} - \frac{\log ( \beta) }{2} \right)  - \frac{X_Q}{2} + (\beta-1) ( \lambda_1 - 2) \notag\\
& - \frac{1}{12} \log(N) - \frac{1}{2} \log ( N^{1/3} ( \beta-1) +1 )+ \O (1).
\end{align}
This yields the claim. \qed

\subsection{Proof of Theorem \ref{thm:main}}
\label{sec:main}

The fact that the quantity,
\beq
Y_N := \frac{N \left(F_N - f(\beta) \right) + \frac{1}{12} \log(N)}{ \sqrt{ \frac{1}{6} \log(N) } }
\eeq
is tight for $\beta$ of the form $\beta = 1 + \alpha \sqrt{ \log(N) } N^{-1/3}$ with $\alpha$ fixed follows immediately from Propositions \ref{prop:fe-ht} and \ref{prop:fe-lt} as well as the fact that the random variables $X_Q / \sqrt{ \log(N) }$ and $N^{2/3} ( \lambda_1 -2 )$ are themselves tight due to Theorem \ref{thm:pq} and the convergence of the latter to the Tracy-Widom$_1$ distribution.

For the Gaussian fluctuations in the case that $\alpha \to 0$, we see first that for any $\eps >0$, by Propositions \ref{prop:fe-ht} and \ref{prop:fe-lt} (and the tightness of the random variable $N^{2/3} ( \lambda_1 -2 )$) that there is a $Q >0$ and a $C_\eps >0$ so that, 
\beq
\left| Y_N +\frac{X_Q}{ \sqrt{ 2 \log(N)/3}} \right| \leq C_\eps \left( \frac{1+ \log( \log(N))}{ \sqrt{ \log(N)}} + |\alpha| \right)
\eeq
with probability at least $1-\eps$ for all $N$ large enough.  Hence, for Lipschitz $O$, (with $Y$ as above)
\begin{align}
\limsup_{N \to \infty} \left| \ee[ O (Y_N) ] - \ee[ O ( -X_Q / \sqrt{ 2 \log(N)/3 } ) ] \right| \leq ||O||_\infty \eps.
\end{align}
On the other hand, 
\beq
\lim_{N \to \infty} \ee[ O ( -X_Q / \sqrt{ 2 \log(N)/3 } ) ] = \ee[ O (Z) ]
\eeq
where $Z$ is a standard normal random variable.  Therefore, taking $\eps \to 0$ yields the claim.

From Proposition \ref{prop:ht-b2} we see that if $\beta = 1 + \alpha \sqrt{\log(N) }N^{-1/3}$ with $\alpha <0$ fixed, then for any $\eps >0$ there is a $C_\eps >0$ and $Q >0$ so that,
\beq
\left| Y_N +\frac{X_Q}{ \sqrt{ 2 \log(N)/3}} \right|  \leq C_\eps \left( \frac{ 1 + \log ( \alpha \log(N) )}{ \sqrt{ \log(N) }} \right).
\eeq
The Gaussian fluctuations for $\alpha <0$ then follows from the same argument as in the $\alpha \to 0$ case.

Finally, we consider the case $\alpha \to \infty$ as $N \to \infty$.  From Proposition \ref{prop:fe-lt} we see that for any $\eps >0$ there is a $C_\eps >0$ on which,
\beq
 \left| \frac{ N (F_N ( \beta) - f ( \beta ) ) + \log(N)/12}{N^{1/3} (\beta -1) } - N^{2/3} ( \lambda_1 - 2) \right| \leq C_{\eps} \left( \frac{1}{ \sqrt{ \log(N) }} + \frac{ 1+ \log (1+ \alpha)}{ \alpha } \right)
\eeq
where we used that $X_Q / \sqrt{ \log(N)}$ is tight.  Hence, if $\alpha \to \infty$ as $N \to \infty$, then for Lipschitz $O$,
\begin{align}
\limsup_{N \to \infty} \left| \ee[ O ( \sqrt{ \frac{1}{6} \log(N) } N^{-1/3} ( \beta-1)^{-1} Y_N ) ] - \ee[ O ( N^{2/3} ( \lambda_1 -2 ) ) ] \right| \leq \eps ||O||_\infty.
\end{align}
The claim follows by the convergence of $N^{2/3} ( \lambda_1 -2)$ to the Tracy-Widom$_1$ distribution. \qed


\bibliography{ssk_bib}{}
\bibliographystyle{abbrv}

\end{document}